\documentclass[english,twocolumn]{article}
\usepackage{csquotes}

\usepackage{xcolor}
\usepackage{amsfonts,amssymb,amsmath,amsthm}
\newcommand{\N}{\mathbb{N}}
\newcommand{\R}{\mathbb{R}}

\usepackage{subcaption}
\usepackage{graphicx}
\usepackage{enumitem}
 \usepackage{caption}
 \usepackage{algorithm,algpseudocode}
\makeatletter
\def\plist@algorithm{Alg.\space}
\makeatother
\usepackage[italicComments=true, beginComment=//~]{algpseudocodex}

\newtheorem{thm}{Theorem}
\newtheorem{lem}{Lemma}

\newtheorem{ass}{Assumption}
\newtheorem{remark}{Remark}
\newtheorem{definition}{Definition}

\usepackage[capitalise]{cleveref}
\crefname{ass}{Assumption}{Assumptions}
\crefname{thm}{Theorem}{Theorems}
\crefname{prop}{Proposition}{Propositions}
\crefname{lem}{Lemma}{Lemmas}
\crefname{cor}{Corollary}{Corollaries}
\crefname{algorithm}{Algorithm}{Algs.}
\crefname{figure}{Fig.\!}{Figs.\!}

\usepackage{tikz}
\usetikzlibrary{patterns,arrows, math}
\usetikzlibrary{shapes, arrows}
\tikzstyle{block} = [draw, fill=white, rectangle, minimum height=3em, minimum width=6em]
\tikzstyle{output} = [coordinate]
\tikzstyle{input} = [coordinate]
\usepackage{pgfplots}
\usetikzlibrary{intersections, pgfplots.fillbetween}
\usetikzlibrary{positioning}
\usetikzlibrary{backgrounds}
\pgfplotsset{compat=1.18} 

\usepackage[backend=bibtex,sortcites,defernumbers,maxnames=99,url=false,doi=false]{biblatex}
\AtEveryBibitem{\clearfield{pages}} 
\bibliography{arXiv_version}

\usepackage[bottom=2cm,top=2cm,left=2cm,right=2cm]{geometry}

\begin{document}

  \title{Two-component controller design to safeguard data-driven predictive control \\ 
\vspace*{+0.5em}
  \large{A tutorial to safe learning-based predictive control: exemplified with DeePC and Koopman MPC}}
  
\author{Lea Bold$^*$, Lukas Lanza$^*$, Karl Worthmann\footnote{Optimization-based Control Group, Institute of Mathematics, Technische Universit\"at Ilmenau, Weimarer Stra\ss e 25, 98693 Ilmenau, Germany. 
 E-mail: \{lea.bold, lukas.lanza, karl.worthmann\}@tu-ilmenau.de \\
\textbf{Corresponding author:} Lukas Lanza
}
}

\date{}

\maketitle

\paragraph{Abstract.}
  We design a two-component controller to achieve reference tracking with output constraints~--~exemplified on systems of relative degree two.
  One component is a data-driven or learning-based predictive controller, which uses data samples to learn a model and predict the future behavior of the system. 
    We exemplify this component concisely by data-enabled predictive control (DeePC) and by model predictive control based on extended dynamic mode decomposition (EDMD).
  The second component is a model-free high-gain feedback controller, which ensures satisfaction of the output constraints if that cannot be guaranteed by the predictive controller. This may be the case, for example, if too little data has been collected for learning or no (sufficient) guarantees on the approximation accuracy derived.
  In particular, the reactive/adaptive feedback controller can 
  be used to support the learning process by leading safely through the state space to collect {suitable} data, e.g., to ensure a sufficiently-small fill distance.
  Numerical examples are provided to illustrate the combination of EDMD-based model predictive control and a safeguarding feedback for the set-point transitions including the transition between the set points within prescribed bounds.

\section{Introduction}
Recently, there has been intensive research on data-driven methods in systems and control, see, e.g., the recent work~\cite{martin:allgower:2023}. 
In this note we focus on the control task of output tracking.
For linear time-invariant (LTI) systems, the so-called fundamental lemma~\cite{willems2005note} by Jan C.~Willems and coauthors paved the way for direct data-based predictive control, see, e.g.,~\cite{coulson2019data} (data-enabled predictive control, DeePC) for discrete-time LTI systems and also the recent survey~\cite{faulwasser2023behavioral} for extensions to stochastic LTI descriptor systems --~provided that the collected input-output data is persistently exciting.
The fundamental lemma has been further extended to, e.g., linear parameter-varying~\cite{verhoek2021fundamental}, flat~\cite{alsalti2021data}, and continuous-time systems~\cite{schmitz2024continuous}.
Alternatives, suitable for nonlinear systems,
are machine-learning methods. 
In particular, reinforcement learning (RL) has emerged as a particularly successful family of machine learning methods due to its close relationship to dynamic programming, see, e.g.,~\cite{kiumarsi2017optimal,bertsekas2019reinforcement}.
Despite its undeniable success, one of the main drawbacks
is the comparatively-large data demand. 
Another, so-called indirect~\cite{dorfler2022bridging}, technique using data-driven surrogate models for predictive control is based on the Koopman theory introduced in~\cite{koopman1931hamiltonian,koopman1932dynamical}.
Herein, the nonlinear system is first lifted using so-called observable functions, before the observables are propagated by the linear, but --~in general~-- infinite-dimensional Koopman operator, see, e.g., the recent review article~\cite{brunton2021modern}.
Using extended dynamic mode decomposition (EDMD; \cite{williams:kevrekidis:rowley:2015,klus:nuske:peitz:niemann:clementi:schutte:2020}) a data-driven numerically-tractable approximation is computed. 
EDMD in the Koopman framework has been applied successfully in various examples, see, e.g., the textbook~\cite{mauroy2020koopman} as well as the references therein.
Recently, this technique was successfully extended to control systems, see, e.g., \cite{brunton2016koopman} and~\cite{surana:2016} for linear and bilinear surrogate models, respectively. 
While linear surrogate models are attractive due to their simplicity, see, e.g., \cite{korda2018linear}, the approximation accuracy is limited as recently shown in~\cite{bruder2021advantages,iacob:toth:schoukens:2022}.
A key advantage of Koopman-based methods is the available error analysis in the infinite-data limit~\cite{korda:mezic:2018b}, which was then further elaborated such that nowadays probabilistic finite-data error bounds exists for dynamical systems~\cite{mezic2022numerical,zhang:zuazua:2023} as well as their extension to stochastic and control systems~\cite{nuske2023finite}, see also~\cite{PhilScha24} for an extension to discrete- and continuous-time Markov processes on Polish spaces. 
Based on such error bounds, a Koopman-based (predictive) controller design with end-to-end guarantees can be achieved~\cite{bold2024data,strasser2024koopman}.
However, uniform error bounds typically require invariance of the dictionary, i.e., the space spanned by finitely-many observables. 
While such assumptions may hold for autonomous systems relying on Koopman modes and eigenfunctions~\cite{mezic:2005,mezic:2013,iacob:toth:schoukens:2022}, the respective conditions~\cite{goswami:paley:2021} for control are more restrictive.
Here, kernel EDMD~\cite{klus:nuske:hamzi:2020,philipp2024error} provides a remedy, which even allows to rigorously establish Koopman invariance of the underlying reproducing kernel Hilbert space (RKHS) and, thus, to uniform error bounds~\cite{kohne2024infty}.
Moreover, the respective techniques were recently extended to control-affine systems in~\cite{BoldPhil24}. 

To achieve such good tracking performances, all of the methods mentioned above rely on an offline training phase meaning that data have to be collected in advance. 
In this paper, we address this aspect.
To ensure constraint satisfaction while taking advantage of learning control, 
the field of safe learning has gained importance and several safety frameworks have been proposed, see, e.g.,
the works~\cite{garcia2015comprehensive,hewing2020learning,osborne2021review} for concise overviews. 
To name but a few particular approaches, 
Hamilton-Jacobi equations and respective reachable sets are analyzed in~\cite{chen2018hamilton},
control barrier functions are utilized in~\cite{ames2019control},
Lyapunov stability analysis results are used in~\cite{perkins2002lyapunov} for RL,
and safe learning in MPC is under consideration in~\cite{aswani2013provably}.
Another approach that has attracted attention recently is the construction of so-called predictive safety filters. 
The idea is to test the controls using a model to ensure compliance with the prescribed conditions, cf.~\cite{wabersich2021predictive,wabersich2022predictive} invoking barrier functions.

In this paper we use the idea of two-component controllers, cf.~\cite{skogestad2005multivariable} or~\cite{araki2003two} for a tutorial introduction.
The idea of combining two control strategies has already been applied successfully, e.g., for an event-based PI controller~\cite{sanchez2011two}, combining feedforward open-loop control with feedback, as reported in~\cite{drucker2024experimental} for a mechanical system, 
in~\cite{schmitz2023safe} for continuous-time sampled-data systems in combination with DeePC using Willems et al.'s fundamental lemma,
or in the superposition of RL-based controller and feedback~\cite{johannink2019residual,lanza2024_sampleddata,gottschalk2024reinforcement} to name but a few. 
The two-component controller proposed in this paper is structurally similar to those in~\cite{lanza2024_sampleddata,gottschalk2024reinforcement} and
consists of a learning-based predictive controller on the one hand, and the so-called funnel controller~\cite{IlchRyan02b,BergIlch21} as a safeguarding reactive feedback controller on the other hand.

In the present paper we address the following two aspects of data-driven and learning-based control. 

\noindent 
\textbf{Achieve a prescribed fill distance.}
We apply the feedback controller with prescribed performance to the system and thereby make the system follow a pre-defined
reference trajectory through the state space. This reference is chosen in such a way that the system visits 
desired pre-specified configurations, i.e., points in the state space, where samples are taken.
If the reference is chosen such that the \textit{fill distance} see \Cref{Def:FillDistance}) of system samples satisfies a predefined threshold, then, following the ideas in~\cite{kohne2024infty,BoldPhil24} to approximate the Koopman operator using kernel EDMD, $L^\infty$-error bounds on the approximation error can be obtained.

\noindent 
\textbf{Online learning predictive control.}
Safeguarded by the reactive feedback controller, a tracking task is started, and system data is collected during runtime, cf.~\Cref{Fig:Controller}.
At some point, e.g., if a certain threshold for the fill distance is ensured or a certain amount of data is available, an EDMD surrogate model is computed and an EDMD-based predictive controller is activated. 
At this stage the reactive feedback controller intervenes only, if the system enters a safety critical region, cf.~\Cref{Fig:SafeAndSafetycritical}.
Data collection is continued (and, potentially, the surrogate model is updated to improve its prediction capability) until, e.g. based on the fill distance, guarantees can be given using solely the predictive controller. 

This paper is organized as follows. \Cref{subsec:system-class:control-objective} introduces the system class under consideration and specifies the control objective.
\Cref{Sec:TwoComponentController} presents our idea of designing a two-component controller; there, we describe the combination of \textit{any} data-driven learning-based controller and the safeguarding outer-loop feedback controller, see~\Cref{Fig:Controller}.
In \Cref{Sec:FunnelControl} we briefly recall the concept of funnel control which is used to safeguard the overall controller.
As a particular instance of a data-driven scheme, we discuss data-enabled predictive control (DeePC) in \Cref{Sec:Willems} and learning of the Koopman operator via Extended Dynamic Mode Decomposition (EDMD) in \Cref{Sec:KoopmanBasedControl}. 
To illustrate the proposed controller, we provide numerical simulations in \Cref{Sec:Numerics}.

\ \\
\textbf{Notation}:
$[n:m] = \mathbb{N}_0 \cap [n,m]$ 
for $n,m \in \N$ with $n \leq m$.
$\langle \cdot, \cdot \rangle$ is the standard inner product on $\R^n$, $\|x\|:=\sqrt{\langle x,x\rangle}$ for $x\in\R^n$; 
and for a symmetric positive definite matrix~$Q \in \R^{n \times n}$ we write $\| x \|_Q := x^\top Q x$.
For an interval $I\subset\R$, 
$C^p(I,\R^n)$ is the set of~$p \in \N$ times continuously differentiable functions~$f: I \to \R^n$,
$L^\infty(I,\R^n)$ is the space of measurable essentially bounded
functions $f: I\to\R^n$ with norm $\|f\|_\infty=\text{esssup}_{t\in I}\|f(t)\|$. 
$W^{k,\infty}(I,\R^n)$ is the Sobolev space of all $k$-times weakly differentiable functions
$f:I\to\R^n$ {with} $f,\dots, f^{(k)}\in L^{\infty}(I,\R^n)$,
and $f_{J}$ denotes the restriction of~$f$ to~$J \subset I$.

\section{Two-component controller design for safe learning} \label{sec:two-component-controller}
In \Cref{subsec:system-class:control-objective}, 
we introduce the system class and 
define the control objective of output tracking within prescribed bounds. 
Then, in \Cref{Sec:TwoComponentController}, we present the main idea, i.e., the design of the two-component controller, which enhances the performance of a safety-ensuring controller by combining it with a learning scheme.
In this paper, we focus on funnel control to reliably ensure tracking within prescribed bounds and, then, use DeePC and EDMD to set up a surrogate model for predictive control to improve the closed-loop behavior.

\subsection{System class and control objective} \label{subsec:system-class:control-objective}
We assume that the system has the following form, which represents many mechanical systems without kinematic constraints, see, e.g.~\cite{bastos2017analysis},
\begin{equation}\label{eq:SystemStructure}
    \begin{aligned}
           \dot x_1(t) &= x_2(t), \\
           \dot x_2(t) &= g_0(x_1(t),x_2(t)) + \sum_{i=1}^{m} g_i(x_1(t),x_2(t)) u_i(t), \\
        y(t) & = h(x_1(t),x_2(t)) = x_1(t),
    \end{aligned}
\end{equation}
where~$y(t) \in \R^m$ is the system's measured output and $u(t) = (u_1(t),\ldots,u_m(t))^\top \in \mathbb{R}^m$ denotes the control at time $t \in \mathbb{R}_{\geq 0}$; 
accordingly, $x(t) \in \R^{2m}$ is the state.
The function~$g_0 \in \mathcal{C}^1(\R^{2m}, \R^{m})$ denotes the drift, and the maps $g_i \in \mathcal{C}^1(\R^{2m},\R^{m})$, $i \in [1:m]$, distribute the input.
For an input~$u \in L^\infty_{\operatorname{loc}}([0,\infty),\mathbb{R}^m)$ the (Carathéodory) solution $x(\cdot;\hat{x},u)$ w.r.t.\ the initial condition $x(0;\hat{x},u) = \hat{x}$ uniquely exists on its maximal interval of existence.
We propose the following structural property on the system~\eqref{eq:SystemStructure}.
\begin{ass}[Known control direction] \label{Ass:G_PosDef}
    The matrix-valued input distribution $G : \R^m \times \R^m \to \R^{m \times m}$ given by $(x_1,x_2) \mapsto [g_1(x_1,x_2), \ldots, g_m(x_1,x_2) ]$ is sign definite.    
\end{ass}
In the following, we assume positive definiteness, i.e., $\langle z, G(x) z \rangle > 0$ for all $z \in \R^m \setminus \{0\}$ and all~$x \in \R^{2m}$.
This is not restrictive since negative definiteness will only change the sign in the feedback law~\eqref{eq:FunnelControl}. 

For the output map~$h(x) = x_1$, system~\eqref{eq:SystemStructure} can be written as input-output dynamics
\begin{align}
    \ddot y(t) &= g_0(y(t),\dot y (t)) + \sum_{i=1}^{m} g_i(y(t),\dot y(t)) u_i(t) \nonumber \\
    & = g_0(y(t),\dot y (t)) + G(y(t), \dot y(t)) u(t), \label{eq:InputOutputSystem}
\end{align}
using~$G$ from \Cref{Ass:G_PosDef}.

\begin{remark}
While systems~\eqref{eq:InputOutputSystem} are of order two,
the proposed method can be extended to more general input-output systems of order~$r \in \N$ of the form    
\begin{equation*}
    \begin{aligned}
        y^{(r)}(t) &= f(t,\textbf{T}(y,\dot y,\ldots,y^{(r-1)})(t), u(t)), \\
        y|_{[-b,0]} &= y^0 \in C^{r-1}([-b,0],\R^m),
    \end{aligned}
\end{equation*}
where~$f \in C(\R_{\ge 0} \times \R^{\eta} \times \R^m, \R^m)$ satisfies a so-called high-gain property, cf.~\cite[{Def.~1.2}]{BergIlch21}, and the operator $\textbf{T} : C(\R_{\ge 0},\R^{rm}) \to L_{\rm loc}(\R_{\ge 0},\R^\eta)$ satisfying~\cite[{Def.~1.1}]{BergIlch21} models bounded though unobservable internal dynamics, hysteresis effects and delays.
However, for the sake of clarity,
we restrict the presentation to systems~\eqref{eq:SystemStructure}, respectively~\eqref{eq:InputOutputSystem}.
\hfill $\diamond$
\end{remark}

\begin{figure}
\begin{subfigure}[t]{0.47\linewidth}
\begin{tikzpicture}[>=stealth,scale=0.51]
\draw[->] (0,0) --node[below]{$t$} (6.3,0);
\draw[->] (0,0) -- node[rotate=90,left=1em,above=0.1]{$\|e\|$} (0,3) ;
\draw[scale=2, dashed, domain=0:3, samples=100 , variable=\x, black] plot ({\x}, {1.2*exp(-\x)+0.05});
\draw[scale=2, domain=0:3, samples=100 , variable=\x, black] plot ({\x}, {abs(cos(deg(0.5*\x))*cos(deg(3*\x))*exp(-0.7*\x)});
\draw (1.4,1.4) -- node[above=0.7em,right=0.1em]{\footnotesize{$1/\sigma$}} (1.7,1.8);
\end{tikzpicture}
\subcaption{Evolution of the norm of the error~$e$ within the boundary~$1/\sigma$.}
\label{Fig:ErrorInFunnel}
    \end{subfigure}
    \quad
    \begin{subfigure}[t]{0.47\linewidth}
    \begin{tikzpicture}[>=stealth,scale=0.51]
\draw[, name path= axis,->] (0,0) --node[below]{$t$} (6.3,0);
\draw[->] (0,0) -- node[rotate=90,left=1em, above=0.3em]{$\|e_2\|/\sigma$} (0,3) ;
\draw[name path=upper, scale=2, dashed, domain=0:3, samples=100 , variable=\x, black] plot ({\x}, {1.2*exp(-\x)+0.05});
\draw[name path=lower, scale=2, dashed, domain=0:3, samples=100 , variable=\x, black] plot ({\x}, {0.75*(1.2*exp(-\x)+0.05});
\tikzfillbetween[of=upper and lower]{red, opacity=0.1};
\tikzfillbetween[of=lower and axis]{green, opacity=0.1};
\draw[scale=2, domain=0:3, samples=100 , variable=\x, black] plot ({\x}, {abs(1.9*sin(deg(1*\x))*cos(deg(0.8*\x))*exp(-0.9*\x)});
\draw (1.4,1.4) -- node[above=0.7em,right=0.1em]{\footnotesize{$1/\sigma$}} (1.7,1.8);
\end{tikzpicture}
\subcaption{Safe (green) and safety-critical (red) region for~$e_2$.}
\label{Fig:SafeAndSafetycritical}
\end{subfigure}
\caption{Schematic illustration of tracking error, funnel boundary as well as safe and safety-critical areas~\cite{gottschalk2024reinforcement}.}
\label{Fig:ErrorAndFunnel}
\end{figure}
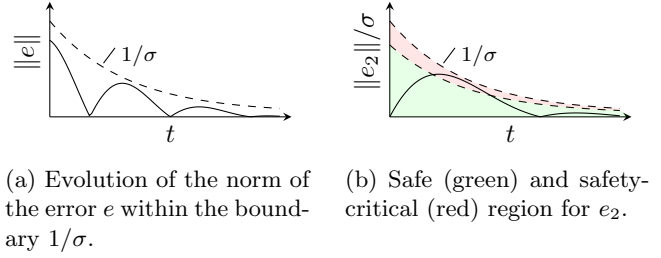

The control objective is that the output~$y$ of system~\eqref{eq:SystemStructure} follows a given reference trajectory $y_{\rm ref} \in W^{2,\infty}([0,\infty),\mathbb{R}^m)$ within predefined error margins for all times. 
The latter means that the tracking error $e(t) := y(t) - y_{\rm ref}(t)$ satisfies
\begin{equation}\label{eq:ControlObjective}
    \|y(t) - y_{\rm ref}(t) \| < 1/\sigma(t) \qquad\forall\, t \ge 0,
\end{equation}
where $\sigma$ is a time-varying error margin
chosen by the control engineer, which belongs to the set
\begin{equation*}
    \Sigma := \{\sigma \in W^{1,\infty}([0,\infty),\R) | \inf_{s \ge 0} \sigma(s) > 0 \}.
\end{equation*}
The control objective~\eqref{eq:ControlObjective} is depicted in \Cref{Fig:ErrorInFunnel}.

\subsection{Two-component controller} \label{Sec:TwoComponentController}
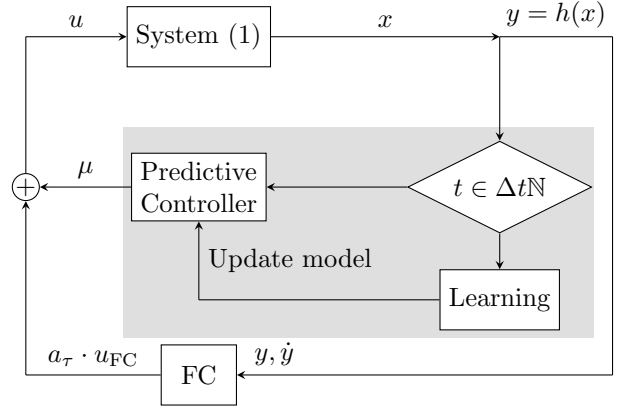
\begin{figure}[htb]
  \begin{center}
\begin{tikzpicture} [auto, node distance=2cm,>=stealth, every text node part/.style={align=center}]
\def\hoch{0.8cm};
\def\breit{1cm};
\def\distu{1.3cm};
\def\dista{1.3cm};
\node [block, minimum width = \breit, minimum height = \hoch,] (System) {System~\eqref{eq:SystemStructure}};
\node [block, minimum width = \breit, minimum height = \hoch, below of=System, node distance = 4.5cm] (FC) {FC};
\coordinate[right of=System, node distance = 4.0cm] (y_mpc) {};
\coordinate[right of=y_mpc, node distance = 1.5cm] (y_fc) {};
\coordinate[below of=y_fc, node distance = 4.5cm] (y_fc2) {};
\node[input, left of=System, node distance = 2.3cm ] (u) {};
\coordinate[right of=u, node distance = 0.0cm] (uin) {};
\node [draw, diamond, aspect=2, fill=white, below of=y_mpc, node distance = 2.0cm] (Sample) {$t \in \Delta t \N$};
\node [block, minimum width = \breit, minimum height = \hoch, below of=Sample, node distance = 1.5cm] (Learn) {Learning};
\node [block, minimum width = \breit, minimum height = \hoch, left of=Sample, node distance = 4.0cm] (MPC) {Predictive \\ Controller};
\coordinate[below of=MPC, node distance = 1.5cm] (Update_MPC) ;
\node[circle,draw=black, fill=white, inner sep=0pt, minimum size=3pt, left of=MPC, node distance = 2.3cm] (uin_control) {$+$} ;
\coordinate[left of=FC, node distance = 2.3cm] (uin_control_FC) ;
\draw[-] (Learn) -- (Update_MPC);
\draw[->] (Update_MPC) -- node[right] {Update model} (MPC);
\draw[->] (u) --  node[above]{$u$} (System);
\draw[->] (System) --  node[above]{$x$} (y_mpc);
\draw[-] (y_mpc) - -node[above]{$y=h(x)$}  (y_fc);
\draw[-] (y_fc) --  (y_fc2);
\draw[->] (y_fc2) --node[above, xshift=-2cm]{$y, \dot y$}  (FC);
\draw[->] (y_mpc) -- (Sample);
\draw[-] (FC) --node [above]{$a_\tau \cdot u_{\rm FC}$} (uin_control_FC)  ;
\draw[->] (uin_control_FC) -- (uin_control)  ;
\draw[->] (MPC) -- node[above]{$\mu$} (uin_control);
\draw[-] (uin_control) -- (uin);
\draw[->] (Sample) -- (Learn);
\draw[->] (Sample) -- (MPC);
\begin{pgfonlayer}{background}
        \fill[gray!25] (5.25,-4) rectangle (-1,-1.2);
    \end{pgfonlayer}
\end{tikzpicture}
\end{center}
 \vspace*{-2mm}
 \caption{Structure of the two-component controller. 
 The data-driven learning-based \textit{Predictive Controller} in the gray box is safeguarded by the \textit{Funnel Controller}. To limit
 activity of the latter, the input~$u_{\rm FC}$ is multiplied by the activation function~$a_\tau$.
 For brevity, the reference~$y_{\rm ref}$ is not shown explicitly but it is provided to both controller components.
 }
 \label{Fig:Controller}
 \end{figure}

To achieve the control objective~\eqref{eq:ControlObjective},
we propose a two-component controller combining a data-driven learning-based controller and a safeguarding feedback controller, cf.~\Cref{Fig:Controller}.
To improve performance, the data-driven controller is typically assumed to be of predictive type (e.g. MPC, DeePC, Subspace Predictive Control).
However, a combination with other reactive feedback controllers, e.g. learning the gain from data~\cite{de2019formulas,de2023learning}, is also possible.

We denote the control signal from the data-driven (predictive) controller by~$\mu$, and the safeguarding control signal by~$u_{\rm FC}$. 
The combination of both control signals~$\mu$ and~$u_{\rm FC}$ is used to achieve the control objective~\eqref{eq:ControlObjective}.
A straight forward combination would be
\begin{equation*}
    u = \mu(t,x) + u_{\rm FC}(t).
\end{equation*}
With this, however, the safeguarding outer-loop funnel controller would interrupt the predictive controller whenever an auxiliary error variable is non-zero.
Hence, we introduce an activation function similar to event-triggered control, see, e.g.,~\cite{lunze2010state,quevedo2014stochastic,liu2015event}.
Following the ideas in~\cite{gottschalk2024reinforcement}, we introduce for a dwell-time~$\tau > 0$ and an activation threshold~$\lambda$ an activation function~$a_\tau(t,e(t),\dot e(t))$, see~\eqref{eq:ActivationFunction} for details.
The activation threshold separates the ``safe'' and ``safety-critical'' region in the output space (more precise in the tracking-error space), cf.~\Cref{Fig:SafeAndSafetycritical}.
Only if the auxiliary signal leaves the safe $\lambda$-region, the funnel controller contributes to the control signal.
Since the funnel controller may instantaneously push the auxiliary variable back into the safe region, the dwell-time~$\tau > 0$ is incorporated to avoid chattering.
For an activation function~$a_\tau(t,e(t),\dot e(t))$, we define the overall control signal 
\begin{equation} \label{eq:ControlSignal}
    u(t) = \mu(t,x) +  a_\tau(t,e(t),\dot e(t)) \, u_{\rm FC}(t).
\end{equation}
To give an idea, how the controller works, we follow the signals in \Cref{Fig:Controller}, omitting the time argument.
Note that (as usual) we design the controller assuming that the prediction step in the data-driven controller (e.g. finite horizon optimal control problem in MPC) is solved instantaneously.
This means that at~$t=k \Delta t$, where~$\Delta t > 0$ is the sampling time, the current state~$x(k\Delta t)$ is given to the prediction scheme which calculates the control input $\mu: [k \Delta t, (k+N) \Delta t) \to \R^m$. Here, $k \in \N$ is the sampling instance, and~$N \in \N$ is the prediction horizon, cf. \Cref{Alg:MPC}. 
We may assume that for the learning scheme the system state~$x$ is available.
If $t = k \Delta t $, 
then the data $x(k \Delta t)$ is given to the learning entity.
Given the data $x(0),\ldots,x(d \Delta t)$ (for some amount of data~$d > 0$) and the recorded input~$u$, the learning entity supplies, e.g., a model that explains the data best (w.r.t certain criteria).
This model is handed over to the predictive controller.
Note that the model is not necessarily updated every~$\Delta t$.
The updates can (and will) be made after arbitrary sampling instances.
Based on the (updated) model the predictive controller produces (e.g. via finite horizon optimal control) a control signal~$\mu$, which is handed over to an adder.
For the funnel controller we take~$y,\dot y$ from the system, and $y_{\rm ref}, \dot y_{\rm ref}$ is given.
These signals are handed over to the funnel controller computing~$u_{\rm FC}$. 
Based on the current error and its derivative, the control signal~$u_{\rm FC}$ is active or inactive by multiplication with~$a_\tau$.
Therefore, the resulting control input is $u = \mu + a_\tau u_{\rm FC}$.
Note that since the funnel controller receives the signals~$y,\dot y$ continuously, the signal~$u_{\rm FC}$ implicitly depends on the signal~$\mu(t,x)$ through the closed-loop continuous-time dynamics~\eqref{eq:SystemStructure},~\eqref{eq:ControlSignal}.

\section{Safeguarding mechanism} \label{Sec:FunnelControl}
As the safeguarding controller component we utilize a high-gain adaptive feedback law.
In the present paper we use \textit{funnel control}, first proposed in~\cite{IlchRyan02b}, but recognize that for the system class under consideration we could have used the \textit{prescribed performance control} methodology~\cite{bechlioulis2010prescribed,bechlioulis2014low} as well.
Both controllers are suitable to achieve~\eqref{eq:ControlObjective} for systems~\eqref{eq:SystemStructure}.
Let a prescribed tolerance~$\sigma \in \Sigma$ be given.
Then, for system~\eqref{eq:SystemStructure} satisfying \Cref{Ass:G_PosDef}, a funnel control feedback law is given in the simple form, cf.~\cite{BergIlch21},
\begin{equation} \label{eq:FunnelControl}
    \begin{aligned}
        e_1(t,y(t)) & : = \sigma(t) e(t) = \sigma(t) ( y(t) - y_{\rm ref}(t) ), \\
        e_2(t,y(t),\dot y(t)) & :=  \sigma(t) \dot e(t) + \frac{e_1(t)}{1-\|e_1(t)\|^2}, \\
        u_{\rm FC}(t,y(t),\dot y(t)) & := - \frac{e_2(t)}{1-\|e_2(t)\|^2},
    \end{aligned}
\end{equation}
where we omit the precise definition of the domain of~$e_1,e_2$, respectively. 
To avoid confusion with the discrete-time setting for EDMD presented in \Cref{Sec:EDMD}, we explicitly state the following assumption.
\begin{ass}[Availability of signals] \label{Ass:SignalsAvailable}
    The signals $y(t)$, $\dot y(t)$, $y_{\rm ref}(t)$ and~$\dot y_{\rm ref}(t)$ are continuously available to the controller~\eqref{eq:FunnelControl}.
\end{ass}

We emphasize the particular feature of funnel control that no system parameters are involved in the feedback law~\eqref{eq:FunnelControl}, i.e., it is model-free. 

\ \\
Next, we make the notion of the activation function mentioned in \Cref{Sec:TwoComponentController} precise.
Let~$\lambda \in (0,1)$ and use the error variables defined in~\eqref{eq:FunnelControl}.
We specify the activation function~$a_\tau : \R_{\ge 0} \times \R^m \to [0,1)$ as follows
\begin{equation} \label{eq:ActivationFunction}
    a_\tau(t,e_2(t)) = \max\{0,  \max_{s \in [t-\tau,t]} \| e_2(s) \| - \lambda \}.
\end{equation}
The activation function has two parameters~$\lambda \in (0,1)$ and~$\tau > 0.$
The activation threshold~$\lambda$ determines the safe region within the funnel boundaries, cf.~\Cref{Fig:SafeAndSafetycritical}. For small values of~$\lambda$ the funnel controller intervenes already for small deviations from the reference. 
For~$\lambda$ close to~$1$ the funnel controller reacts late which results in high peaks of the control input since the feedback has to push back the error rapidly.
The latter is, in particular, relevant when performing numerical simulations.
The parameter~$\tau > 0$ defines the dwell-time to avoid chattering when~$e_2$ leaves the~$\lambda$-region. 
Large values of~$\tau > 0$ keep the funnel controller active for a long time even if~$e_2$ is already in the safe region which may have effects on the prediction based controller.
The more reliable the predictive controller is, the smaller the dwell-time~$\tau$ can be chosen.
Note that, since~$e_2 = e_2(t,y,\dot y)$, the above definition of~$a_\tau$ is in line with the previous informal introduction, where we used $a_\tau(t,e(t),\dot e(t))$. To save on notation, we use~$a_\tau(t,e_2)$ in the following.
As a first result we record the following safeguarding property. 
It combines --~with minor changes~-- the idea of using an activation function as in~\cite{berger2024robust} with the observations shown in~\cite[Thm.~1]{gottschalk2024reinforcement}.
The feasibility proof mainly relies on the proof of~\cite[Thm.~1.9]{BergIlch21}, and the insights in~\cite[Thm.~5.1]{lanza2024_sampleddata} for sampled-data systems.
\begin{thm}[Safeguarding property] \label{Lem:OutputConstraintSatisfaction}
    Consider a system~\eqref{eq:SystemStructure}. Let a reference trajectory $y_{\rm ref} \in W^{2,\infty}(\R_{\ge 0},\R^m)$ and a funnel function~$\sigma \in \Sigma$ be given.
    If \Cref{Ass:G_PosDef,Ass:SignalsAvailable} are satisfied and for the auxiliary variables in~\eqref{eq:FunnelControl} it holds
    \begin{equation*}
        \|e_1(0)\| < 1 \ \text{and} \ \|e_2(0)\| < 1,
    \end{equation*}
    then any solution of the closed loop system~\eqref{eq:SystemStructure},~\eqref{eq:ControlSignal} satisfies $\|y(t) - y_{\rm ref}(t) \| < 1/\sigma(t)$ for all~$t \ge 0$, i.e., the control objective~\eqref{eq:ControlObjective} is satisfied, where~$\mu(\cdot)$ is a bounded input generated by the predictive controller, and~$a_\tau(\cdot) \,   u_{\rm FC}(\cdot)$ is given by~\eqref{eq:FunnelControl},~\eqref{eq:ActivationFunction}.
\end{thm}
\textit{Main ideas of the proof.}
    We omit the technical details here and justify the result by following the reasoning given in~\cite[Thm.~1]{gottschalk2024reinforcement} instead.
    The input signal~$\mu(\cdot)$ generated by the prediction-based controller is bounded. 
    For~$\lambda \in (0,1)$ from~\eqref{eq:ActivationFunction} and invoking \Cref{Ass:G_PosDef,Ass:SignalsAvailable} it can be shown by contradiction that there exists $\tilde \varepsilon \in [\lambda,1)$ such that $\|e_i(t)\| \le \tilde \varepsilon$ for all~$t \ge 0$, and $i=1,2$.
    The respective analysis can be restricted to time intervals, where $\|e_2(t)\| \in [\lambda,1)$ (red area in~\Cref{Fig:SafeAndSafetycritical}).
    In these intervals, the analysis in the proof of \cite[Thm.~1.9]{BergIlch21} applies. Therefore, multiplying~$u_{\rm FC}$ by~$\alpha_\tau$, which is potentially zero, does not jeopardize the analysis in the safety-critical zone, meaning that we can employ the same line of arguments as in~\cite{BergIlch21,lanza2024_sampleddata,gottschalk2024reinforcement}. \hfill $\diamond$

\ \\
Note that \Cref{Lem:OutputConstraintSatisfaction} means that the proposed two-component controller design can be used for safe learning during runtime, more precisely, to perform safe learning-based predictive control without offline training, since constraints satisfaction in the output is guaranteed by the safeguarding feedback controller.
This is illustrated by a numerical simulation in \Cref{Sec:Numerics}.

\section{Data-based control}
In this section, we present two instances of a two-component controller.
In \Cref{Sec:Willems}, we recap Willems and coauthors' fundamental lemma~\cite{willems2005note} and its use for predictive control of linear time-invariant systems following~\cite{coulson2019data}. Hereby, we leverage our previous works~\cite{schmitz2023safe,lanza2024_sampleddata}.
Then, we present recent results for nonlinear systems using the Koopman operator in \Cref{Sec:KoopmanBasedControl}.

\subsection{Data-enabled predictive control} \label{Sec:Willems}
Consider the discrete-time LTI system
\begin{equation} \label{eq:DiscreteLTI}
    \begin{aligned}
        x^+ &= Ax + B u, \quad 
        y = C x,
    \end{aligned}
\end{equation}
with output~$y \in \R^m$, where according to system class~\eqref{eq:SystemStructure} and \Cref{Ass:G_PosDef}
\begin{equation*}
    A = \begin{bmatrix}
        0&I\\A_1&A_2
    \end{bmatrix}, 
    \quad B = \begin{bmatrix}
    0 \\B_1    
    \end{bmatrix}, \quad 
    C = \begin{bmatrix}
        I_m & 0
    \end{bmatrix} 
\end{equation*}
for~$A_1,A_2 \in \R^{m \times m}$, and positive definite~$B_1 \in \R^{m \times m}$.
Therefore, system~\eqref{eq:DiscreteLTI} is minimal, i.e., controllable and observable.
To recap DeePC~\cite{coulson2019data}, we first recall the concept of persistently exciting signals and a statement for systems~\eqref{eq:DiscreteLTI}.
For~$1 \le d \in \N$ a real sequence $(u_k)_{k=0}^{d-1}$ is called persistently exciting of order~$L\in \N$ if the so-called Hankel matrix 
\begin{equation*}
    \mathcal H_L(u_{[0,d-1]}) = 
        \begin{bmatrix}
        u_0 & u_1 & \dots & u_{d-L}\\
        u_1 & u_2 & \dots & u_{d-L+1}\\
        \vdots & \vdots & \ddots &\vdots\\
        u_{L-1} & u_L & \dots & u_{d-1}
    \end{bmatrix}
\end{equation*}
has full row rank.
Given a persistently exciting input, the following result can be shown for systems~\eqref{eq:DiscreteLTI} known as the fundamental lemma~\cite{willems2005note}; slightly simplified for the current purpose.
\begin{lem}[\cite{willems2005note}]
 \label{lem:fl}
    Let $((\hat u_k)_{k=0}^{d-1}, (\hat y_k)_{k=0}^{d-1})$ be input-output data of~\eqref{eq:DiscreteLTI}, where the input $(\hat u_k)_{k=0}^{d-1}$ is persistently exciting of order $L+2m$. 
    Then $((u_k)_{k=0}^{L-1}, (y_k)_{k=0}^{L-1})$ is an input-output trajectory of \eqref{eq:DiscreteLTI} of length~$L$ if and only if there exists $\nu\in \R^{d-L+1}$ such that
    \begin{equation*} 
        \begin{bmatrix}
            u_{[0,L-1]}\\
            y_{[0,L-1]}
        \end{bmatrix} =\begin{bmatrix}
            \mathcal H_L(\hat u_{[0,d-1]})\\
            \mathcal H_L(\hat y_{[0,d-1]})
        \end{bmatrix} \nu.
    \end{equation*}
\end{lem}
This representation can be used to set up a purely data-driven predictive controller to make the output of system~\eqref{eq:DiscreteLTI} track a reference trajectory given by $(y_{\text{ref}}(k))_{k=0}^\infty$, while~--~like in MPC~--~satisfying input constraints.
Assume that for a prediction horizon~$N > 0$ the input-output data $((\hat u_k)_{k=0}^{d-1}, (\hat y_k)_{k=0}^{d-1})$ of system~\eqref{eq:DiscreteLTI} is given and that the input data $(\hat u_k)_{k=0}^{d-1}$ is persistently exciting of order $N+4m$. 
In every discrete time step~$k \in \N$ the optimal control problem to be solved is given by (here and later we use quadratic stage costs costs to keep the presentation simple)
\begin{subequations}
\label{eq:DeePC}
\begin{equation}
\label{eq:ocp1}
    \operatorname*{minimize}_{(u,y,\nu)} \sum_{i=k+1}^{k+N} \Bigl(\| y(i)-y_{\text{ref}}(i) \|_Q^2 + \| u(i) \|_R^2\Bigr) 
\end{equation}
with $(u,y) = ((u(i))_{i=k-2m+1}^{k+N},(y(i))_{i=k-2m+1}^{k+N})$ subject to
\begin{align}
\label{eq:ocp2}
    \begin{bmatrix}
        u_{[k-2m+1,k+N]}\\
        y_{[k-2m+1,k+N]}
    \end{bmatrix}& = \begin{bmatrix}
        H_{N+2m}(\hat u)\\H_{N+2m}(\hat y)
    \end{bmatrix}\nu,\\
    \label{eq:ocp3}
    \begin{bmatrix}
        u_{[k-2m+1,k]}\\
        y_{[k-2m+1,k]}
    \end{bmatrix} &= \begin{bmatrix}
        \tilde u_{[k-2m+1,k]}\\
        \tilde y_{[k-2m+1,k]}
    \end{bmatrix},\\
    \label{eq:ocp5}
    \| u(i)\| &\leq u_\text{max},\ i=k+1,\dots, k+N. \nonumber
\end{align}
\end{subequations}
Here, $Q,R\in\R^{m\times m}$ in~\eqref{eq:ocp1} are symmetric positive definite weighting matrices.
The past input-output trajectory is denoted by $(\tilde u,\tilde y)=((\tilde u(i))_{i=k-2m+1}^{k}, (\tilde y(i))_{i=k-2m+1}^{k})$, and in combination with the observability of~\eqref{eq:DiscreteLTI} the constraint~\eqref{eq:ocp3} serves as initial condition such that the trajectory is continued. 
  The main difference to model predictive control is that in the optimal control problem~\eqref{eq:DeePC} the state-space model~\eqref{eq:DiscreteLTI} is replaced by the non-parametric description~\eqref{eq:ocp2} based on \Cref{lem:fl}. \\
For the sake of clear presentation let us denote the control generated by~\eqref{eq:DeePC} by $u_{\rm DPC}$, and $u_{\rm DPC}(t) = u_{\rm DPC}(k)$ for $t \in [k,k+1)$.
According to the idea presented in \Cref{Sec:TwoComponentController} the overall input reads
\begin{equation*}
    u(t) = u_{\rm DPC}(t) + a_\tau(t,e_2) u_{\rm FC}(t).
\end{equation*}
Since \Cref{Lem:OutputConstraintSatisfaction} ensures satisfaction of the control objective~\eqref{eq:ControlObjective}, the data required to set up the Hankel matrices~\eqref{eq:ocp2} can be collected during runtime without offline training.
Numerical results for the combination of the data-driven scheme~\eqref{eq:DeePC} with a (sampled-time) funnel controller are presented in~\cite{schmitz2023safe,lanza2024_sampleddata}, where also adaption of the prediction horizon based on the available data is discussed.

\subsection{Koopman-based control} \label{Sec:KoopmanBasedControl}
Although there exist results extending the techniques based on the fundamental lemma to, e.g., linear parameter-varying~\cite{verhoek2021fundamental}, flat~\cite{alsalti2021data}, and continuous-time systems~\cite{schmitz2024continuous}, in its roots this technique is restricted to linear systems.
Now, due to its ability to handle nonlinear systems, 
we utilize MPC based on an EDMD surrogate model in the prediction step, see, e.g., \cite{korda2018linear,bold2024data} as a second particular instance for the learning-based predictive control component (gray box in \Cref{Fig:Controller}). 
Thereby, we 
take advantage of recently-derived finite-data pointwise error bounds for control systems~\cite{BoldPhil24}.

\subsubsection{Extended Dynamic Mode Decomposition} \label{Sec:EDMD}
First, we briefly recap the basic idea of extended dynamic mode decomposition (EDMD). 
To this end, we consider a discrete-time system
\begin{align}\label{eq:time discr sys_nocontrol}
    x^+ = F(x)
\end{align}
with continuous nonlinear map $F:\mathbb{R}^n \rightarrow \R^n$. 
Then the Koopman operator is defined by the identity 
\begin{align*}
    (\mathcal{K} \varphi)(\tilde{x}) = \varphi(F(\tilde{x}))
\end{align*}
for all $\tilde{x} \in \mathbb{R}^n$, and functions $\varphi: \R^n \rightarrow \R$, the so-called observables. 
The Koopman operator~$\mathcal{K}$ is a bounded linear operator on a suitably-defined space of observable functions~$\varphi$. However, the Koopman operator is, in general, infinite dimensional.

EDMD is a data-based learning method to approximate the compression~$P_\mathbb{V} \mathcal{K}|_\mathbb{V}$ of the Koopman operator on the compact set $\mathbb{X} \subset \R^n$, see, e.g., \cite{williams:kevrekidis:rowley:2015,nuske2023finite} for details. 
Therein, the space $\mathbb{V} := \operatorname{span}\{\psi_j\}_{j = 0}^M$ is spanned by finitely many observables. 
Let $\Psi = (\psi_0, \psi_1, \dots, \psi_M)^\top$ be the vector of all observables with $\psi_0 \equiv 1$ and $\psi_j \in \mathcal{C}^1(\R^n, \R)$ for all $j \in [1:M]$ such that $\Psi: \R^n \rightarrow \R^{M+1}$. 
Then, for randomly i.i.d.\ sampled data pairs~$(x_i,F(x_i))$, $i \in [1:d]$, the operator $P_\mathbb{V} \mathcal{K}|_\mathbb{V}$ is approximated by a linear map, which we represent by a matrix $K = K_d^M \in \R^{(M+1) \times (M+1)}$, where~$d$ refers to the number of data points. 
To this end, the regression problem 
\begin{equation*}
    K = \operatorname{arg min}_{\hat K \in \R^{(M +1)\times (M+1)}} \| \Psi_{X^+} - \hat K \Psi_X \|_F
\end{equation*}
is solved with data matrices $\Psi_X := [\Psi(x_1) \ \dots \ \Psi(x_d)]$ and $\Psi_{X^+}:= [{\Psi(F(x_1)) \ \dots \ \Psi(F(x_d))}]$. 
Using the pseudo-inverse, the solution reads 
$K = \Psi_{X^+} \Psi_X^\dagger$.
Based on this approximation in the lifted space~$\mathbb{V}$ and assuming full rank of the lift~$\Psi$, an EDMD-based surrogate model of system \eqref{eq:time discr sys_nocontrol} can be defined by $x^+ = P_\mathbb{X} {K}\Psi(x)$,
where $P_\mathbb{X}$ is the orthogonal projection on the state space. To simplify the \textit{projection back} onto the state space~\cite{mauroy2016linear,mauroy2019koopman}, we assume that the coordinate functions are observables, i.e., $\psi_i(x) = e_i^\top x$ for $i \in [1:n]$, where $e_i$ denotes the $i$-th unit vector, and refer to~\cite{goor:mahony:schaller:worthmann:2023} for more-sophisticated reprojection techniques. 
Then, using the simple \textit{coordinate reprojection} yields the surrogate model
\begin{align} 
    x^+ = \left[ \begin{array}{cc}
        I_n & 0_{n \times m}
    \end{array} \right] K\Psi(x). 
\end{align}
\begin{remark}\label{rem:Koopman:continuous-time}
    If the underlying dynamics are governed by a continuous-time dynamical system $\dot{x}(t) = f(x(t))$ with locally Lipschitz-continuous map $f:\mathbb{X} \rightarrow \R^n$, the corresponding discrete-time system is given by
    \begin{align*} 
        x^+ = F(\tilde{x}) := \tilde{x} + \int_0^{\Delta t} f(x(s;\tilde{x}))\,\mathrm{d}s
    \end{align*}
    for a fixed time step $\Delta t > 0$, where $x(\cdot; \tilde{x})$ represents the solution of the (ordinary) differential equation satisfying the initial condition $x(0;\tilde{x}) = \tilde{x}$. 
\end{remark}

\subsubsection{EDMD-based predicted control}\label{subsec:EDMD:control}
Next, we briefly recap an extension of the presented approximation technique EDMD to nonlinear control-affine systems. 
There are different approaches to incorporate inputs, e.g., \textit{EDMD with control} (EDMDc) corresponding to a linear data-driven surrogate model~\cite{proctor2016dynamic,korda2018linear}, 
or the recently proposed method in~\cite{asada2024control} for systems~\eqref{eq:SystemStructure}, where the input enters linearly.
Here, we proceed with a bilinear model, see, e.g., \cite{surana:2016,peitz:otto:rowley:2020} and the references therein. 
The reasoning behind is twofold. On the one hand, the authors of~\cite{iacob:toth:schoukens:2022} clearly point out structural limitations of linear surrogate models due to the imposed decoupling of (lifted) state and control. 
On the other hand, they proposed LPV systems as a remedy pointing out that bilinear models suffice if the dynamics of the original system are control affine. 
To this end, the key insight is that the Koopman generator, which generates the semigroup $(\mathcal{K}^t)_{t \geq 0}$ of linear Koopman operators, preserves control affinity. 
This property is approximately preserved for the Koopman operator for sampled-data systems with zero-order hold, i.e., 
\begin{equation}\nonumber
    x^+ = \tilde{x} + \int_0^{\Delta t} g_0(x(s;\tilde{x},u))\,\mathrm{d}s + \sum_{i = 1}^m \tilde{u}_i \int_0^{\Delta t} g_i(x(s;\tilde{x},u))\,\mathrm{d}s
\end{equation}
with $u(s) \equiv \tilde{u} \in \mathbb{R}^m$, and results in an additional error of magnitude~$\mathcal{O}(\Delta t^2)$, see \cite[Rem.~4.2]{BoldPhil24} for details. 
Overall, using the approximate control affinity analogously to~\cite{peitz:otto:rowley:2020,nuske2023finite}, we obtain the bilinear EDMD surrogate
\begin{align}\label{eq:sur:control}
    x^+ = \widehat{F}(x, u) = \left[ \begin{array}{cc}
        I_n & 0_{n \times m}
    \end{array} \right] {K}_u^{\Delta t} \Psi(x)
\end{align}
to approximate the Koopman operator $\mathcal{K}_u^{\Delta t}$ for $u \in \R^m$ using the decomposition
\begin{align*}
    K_u^{\Delta t} & = K_0^{\Delta t} + \sum_{i=0}^m u_i (K_i^{\Delta t} - K_0^{\Delta t}),
\end{align*}
where $K_0^{\Delta t}$ corresponds to $\tilde{u} = 0$ and $K_i^{\Delta t}$ to $\tilde{u} = e_i$, $i \in [1:m]$. 
Koopman-based surrogate models can be used to perform model predictive control based on system data, as shown, e.g., in \cite{korda2018linear} for EDMDc or with the previously introduced methodology resulting in a bilinear surrogate~\eqref{eq:sur:control}, in \cite{folkestad2021koopman, bold2024data, rosenfelder2024data, narasingam2023data} to name a few. 
The underlying idea is to perform predictions for the observables satisfying the discrete-time dynamics~\eqref{eq:sur:control}. 
In the following, $\bar{x}_{\bar{u}}(\kappa, \hat{x})$ denotes the predicted trajectory emanating from state~$\hat{x}$ at time $\kappa$, $\kappa \in [0:N]$, if the sequence of control values $(\bar{u}(\kappa))_{\kappa = 0}^{N-1}$ is applied. 
To keep the presentation technically simple, we use the quadratic stage costs $\ell: \N_{0} \times \R^n \times \R^m \rightarrow \R_{\geq 0}$ defined by
\[
    \ell(k,x, u) := \|x - x_{\text{ref}}(k \Delta t)\|^2_Q + \| u \|^2_R
\]
with symmetric and positive definite matrices $Q \in \R^{n \times n}$ and $R \in \R^{m \times m}$.
An EDMD-based MPC scheme is presented in \Cref{Alg:MPC}.

\begin{algorithm}[ht!]
\caption{EDMD-based MPC}
\label{Alg:MPC}
\begin{algorithmic}
\Require{Time instant~$k \in \N_0$, 
    prediction horizon~$N \in \mathbb{N}_{\geq 2}$, time shift $\Delta t > 0$}
    \State 1. At time instant~$k \in \mathbb{N}_0$ measure $\hat x = x(k \Delta t)$ and receive the reference $(x_{\rm ref}(\kappa))_{\kappa=k}^{k+N}$
    \State 2. Solve the optimization problem 
    \begin{subequations}\label{eq:OCP}
    \begin{alignat*}{3}%
            &\hspace{0pt} \underset{{\bar{u} =(\bar{u}(\kappa))_{\kappa=0}^{N-1}}}{\textnormal{minimize}}
            &&\hspace{6pt}\!\sum_{\kappa=0}^{N-1}\!\ell(k + \kappa,\bar{x}_{\bar{u}}(\kappa, \hat{x}),{\bar{u}}(\kappa))\\
            & \textnormal{subject to}
            &&\hspace{6pt} \bar{x}_{\bar{u}}(\kappa+1\,, \hat{x}) = \widehat{F} \left(\bar{x}_{\bar{u}}(\kappa\,, \hat{x}), \bar u (\kappa) \right),\\
            &&&\hspace{6pt} \bar{u}(\kappa)\in\mathbb{U},\quad \forall \ \kappa \in [0:N-1],
    \end{alignat*}
    \end{subequations}
    to compute a minimizing sequence $\bar{u}^\star  := (\bar{u}^\star(\kappa))_{\kappa=0}^{N-1}$ of control values
    \State 3. Apply feedback value $\mu_{\text{MPC}}(k,x(k \Delta t)) = \bar{u}^\star(0)$, set $k \gets k+1$, and go to step~1
\end{algorithmic}
\end{algorithm}
It has been shown in~\cite{bold2024data} that under suitable assumptions \Cref{Alg:MPC} defines a practical asymptotically stabilizing controller for nonlinear systems~\eqref{eq:SystemStructure} using finite-data error bounds, cost controllability, and a sufficiently long prediction horizon, see \cite{grune2010analysis} and \cite{boccia2014stability}. 

\begin{remark} \label{rem:aTauZero}
Note that the EDMD procedure relies on piecewise constant control inputs~$u(t) \equiv u_k$ for ${t \in [k \Delta t, (k+1) \Delta t)}$.
Therefore, to obtain reasonable data to learn from, collected during runtime, the data collection can be equipped with a post-processing condition.
One possibility is to apply~$u_k$ for~$t \in [k \Delta t, (k+1) \Delta t)$ and evaluate whether the funnel controller was active in this interval or not by considering the activation function in~\eqref{eq:ActivationFunction}.
If $\alpha_\tau(t,e_2(t)) = 0$ for all~$t \in [k \Delta t, (k+1) \Delta t)$, then the data can be used for learning.
\hfill $\diamond$
\end{remark}

As shown e.g. in \cite{bold2024koopman} a key challenge for successive application of EDMD and EDMD-based MPC is the choice of the observable functions constituting the dictionary that spans the subspace~$\mathbb{V}$ where the Koopman operator is approximated on. 
The selection and construction of dictionaries that suit the underlying system is discussed e.g. in \cite{korda2020optimal}. 
In the next section we concisely recap kernel-based EDMD, see \cite{williams2014kernel}, which attends to the aspect of choosing a dictionary. 
In these works, the dictionary is constructed from kernels, which are defined through the available data points.

\subsubsection{Sampling and approximation error} \label{Sec:PrescribedApproximationError}
Performing any kind of approximate identification or learning of dynamics from data, the question naturally arises as to how good the approximation is.
In recent years this question has been addressed intensively in the context of the Koopman  
framework.
The bounds on the approximation error of the Koopman operator/generator are typically provided as probabilistic statements, see \cite{mezic2022numerical,zhang:zuazua:2023,PhilScha24} and the references therein.
In particular, the error bounds are derived on the assumption of either i.i.d. or ergodic sampling.
However, in practice it is difficult to meet the requirements for i.i.d. or ergodic sampling. 
This issue can be avoided by
using kernel-EDMD to approximate autonomous dynamics using the canonical features as observables in a suitably-chosen reproducing kernel Hilbert space (RKHS).
Moreover, the respective error bounds provide pointwise bounds exploiting the so-called reproducing property and Koopman invariance of the underlying RKHS, see~\cite{kohne2024infty}. 
Then, the results were extended in~\cite{BoldPhil24} to control systems, where, in addition, proportional error bounds w.r.t.\ the 
\textit{fill distance} of data samples were established.
The fill distance~$h_\mathcal{X}$ corresponds to the radius of the largest ball with center in~$\mathbb{X}$, which has no intersection with the set of samples~$\mathcal{X}$.
\begin{definition}[Fill distance~$h_\mathcal{X}$] \label{Def:FillDistance}
    For a bounded set ${\mathbb{X} \subset \mathbb{R}^n}$ and a set of samples $\mathcal{X} := \bigcup_{i=1}^{d} \{ x_i \} \subset \mathbb{X}$, $d \in \mathbb{N}$, the fill distance~$h_{\mathcal{X}}$ w.r.t.~$\mathbb{X}$ is defined by
    \begin{equation*}
        h_{\mathcal{X}} := \sup_{x \in \mathbb{X}} \min_{x_i \in \mathcal{X}} \| x- x_i\|.
    \end{equation*}
\end{definition}
The first uniform bound for the full approximation error for kernel EDMD (kEDMD) was derived in~\cite{kohne2024infty} using the RKHS generated by Wendland kernel functions, see~\cite{wendland2004scattered}. 
We recall the result (in a slightly-simplified fashion) in the following theorem.
\begin{thm}[Theorem~5.2 of \cite{kohne2024infty}] \label{thm:koehne}
    Let the bounded set~$\mathbb{X} \subset \mathbb{R}^n$ have Lipschitz boundary. Further, let $\mathbb{H}$ be the RKHS generated by Wendland kernels with smoothness degree $k \in \N$.
    Further, let $\mathcal{X} := \bigcup_{i=1}^{d} \{ x_i \} \subset \mathbb{X}$ be a finite set of pairwise-distinct samples and ${F \in C^p(\mathbb{X},\R^n)}$, where $p = \lceil \frac{n + 1}{2} + k\rceil $, for the right-hand side of system~\eqref{eq:time discr sys_nocontrol}. 
    Then, there exist constants $C,h_0 > 0$ such that the following bound on the full approximation error holds for all fill distances $h_\mathcal{X} \leq h_0$
    \begin{align*}
        \| \mathcal{K} - K\|_{\mathbb{H} \rightarrow C^p(\mathbb{X},\R^n)} \le Ch_\mathcal{X}^{k+1 \addslash 2}.
    \end{align*}
\end{thm}
In conclusion, \cref{thm:koehne} states that the kEDMD surrogate~$K$ approximates the dynamics (of an autonomous dynamical system) arbitrarily well if only the fill distance~$h_\mathcal{X}$ 
is sufficiently small.
This result was extended in~\cite{BoldPhil24} to control-affine systems 
\begin{equation} \label{eq:DiscreteSystem}
    x^+ = \hat{F}(x,u) := \hat g_0(x) + \sum_{j=1}^m \hat g_j(x) u_j
\end{equation} 
with $\hat g_j \in \mathcal{C}(\R^n,\R^n)$, $j \in [0:m]$.
We briefly recap the main ingredients analogously to~\cite{BoldScha25}. \\
\textbf{Data.} 
Let~$\mathcal{X} = \{x_1,\ldots,x_{d}\} \subset \mathbb{X}$ be a set of pairwise distinct points, the so-called virtual-observation points with \textit{cluster} radius~$\varepsilon_c > 0$.
For each cluster~$i \in [1:d]$, let~$d_i \ge m + 1$ data triplets $(x_{ij},u_{ij},x_{ij}^+) \in B_{\varepsilon_c}(x_i) \times \mathbb{U} \times \mathbb{R}^n$ satisfying $x_{ij}^+ = \hat{F}(x_{ij},u_{ij})$ and the condition $\text{rank}([u_{i1}\, \cdots\, u_{id_i}]) = m$, reminiscent on persistency of excitation as used in~\cite{willems2005note}, on the input data. 
Then, we approximate the function values $\hat{g}_0(x_i),\ldots,\hat{g}_m(x_i)$ in~\eqref{eq:DiscreteSystem}
at the virtual observation points~$x_i \in \mathcal{X}$ by solving the linear regression problem
\begin{equation*}
    [\tilde{g}_0(x_i), \tilde{G}(x_i)] = 
    \underset{H \in \R^{n \times (1+m)}}{\text{argmin}} \left\| [x_{i1}^+ \,  \cdots \, x_{id_i}^+] - H U_i \right\|_F,
\end{equation*}
where
$
    U_i := \begin{bmatrix}
        1 & \cdots & 1 \\
        u_{i1} & \cdots & u_{id_i}
    \end{bmatrix}.
$
Due to the imposed rank condition, the unique solution can be expressed using the pseudo-inverse~$U_i^\dagger$.

\noindent
\textbf{Interpolating the dynamics.} 
The canonical feature maps~$\Phi_{x_i} = \text{k}(x_i,\cdot)$ of the Wendland kernels are used to approximately represent the evolution of an observable~$\psi : \R^n \to \R$ under control input~$u$ by
\begin{equation*}
    \psi(x^+) = \sum_{i=1}^d \big[ (\hat K_0 \psi_{\mathcal{X}})_i \Phi_{x_i}(x) 
    + \sum_{j=1}^m (\hat K_j \psi_{\mathcal{X}})_i \Phi_{x_i}(x) u_j 
    \big],
\end{equation*}
where we use the notation $\psi_\mathcal{X} = (\psi(x_1), \ldots, \psi(x_d))^\top$ and $\hat K_j = K_{\mathcal{X}}^{-1} K_{{\tilde{g}}_j} K_{\mathcal{X}}^{-1} $ with $K_{\mathcal{X}} = (\text{k}(x_i,x_j))_{i,j=1}^d$ 
and $K_{{\tilde{g}}_j} = (\text{k}(x_i, \tilde{g}_j(x_i)))_{i=1}^{d}$ for $j \in [0:m]$.
Choosing the coordinate functions as observables yields a data-driven surrogate for~\eqref{eq:DiscreteSystem} given by 
\begin{equation*}
    x^+ = {F^\varepsilon(x,u) :=} g_0^\varepsilon(x) + \sum_{j=1}^m g_j^\varepsilon (x) u_j,
\end{equation*}
where~$\varepsilon$ refers to the notation in the following result.
\begin{thm}[Theorem~2 of \cite{BoldScha25}] \label{Thm:kEDMDErrorBounds}
    For Wendland kernels of smoothness degree~$k \ge 1 + \tfrac{(-1)^n+1}{2}$, there exist constants~$C,h_0, \bar G, \bar U > 0$ such that for all fill distances $h_{\mathcal{X}} \le h_0$ the approximation error satisfies for all~$(x,u) \in \mathbb{X} \times \mathbb{U}$
\begin{align*}
    \| \hat F(x,u) - F^\varepsilon(x,u) \|_\infty
    &\le 
    C( \bar G h_{\mathcal{X}}^{k-\tfrac{1}{2}} \operatorname{dist}(x,\mathcal{X}) + \bar U \varepsilon_c),
\end{align*}
where~$\bar G$ depends on the approximations~$[\tilde g_0,\ldots, \tilde g_m]$, and~$\bar U$ depends on the input data, the kernel functions and the matrix~$K_{\mathcal{X}}^{-1}$, and~$\varepsilon_c$ is the size of the clusters.
\end{thm}
For a detailed representation of the constants we refer to~\cite[Theorem~2]{BoldScha25}.

\begin{figure}[H]
\centering
\begin{tikzpicture}[node distance=1.4cm,
    every node/.style={fill=white}, align=center, rounded corners]
  \node (FillDistance_r)[block]{
  Given: virtual-observation points~$\mathcal{X} = \{x_i\}_{i=1}^d \subset \mathbb{X}$ \\ 
  \phantom{Given:} with fill distance~$h_{\mathcal{X}}$, cluster radius~$\varepsilon_c > 0$} ;
  \node (FindReference)[block,anchor=north,  below of=FillDistance_r, node distance = 1.6cm]{
  Define reference~$y_{\rm ref}: \R_{\ge 0} \to \R^m$ \\ and sampling time~$\Delta t> 0$ such that \\
  $ \bigcup_{i=1}^d\{ (
    y_{\rm ref}(i \Delta t)^\top, \dot y_{\rm ref}(i \Delta t)^\top)^\top \} = \mathcal{X}$ } ;
  \node (Control) [block,anchor=north, below of=FindReference, node distance = 1.6cm] {
  Apply controller~\eqref{eq:FunnelControl} to system~\eqref{eq:SystemStructure} with $\sigma \ge 3/\varepsilon_c$ \\
  to achieve $x(i \Delta t) \in B_{\varepsilon_c}(x_i)$ in virtue of~\eqref{eq:FunnelBoundariesStatespace}
  };
    \node (Learn) [block,anchor=north, below of=Control]{
    EDMD-based surrogate model satisfies \\ approximation bounds in \Cref{Thm:kEDMDErrorBounds}
    };
  \draw[->] (FillDistance_r) -- (FindReference);
  \draw[->] (FindReference) -- (Control);
  \draw[->] (Control) -- (Learn);
    \end{tikzpicture}
    \caption{Flowchart illustrating the procedure to generate a data-driven surrogate model with prescribed accuracy using~\eqref{eq:FunnelControl}}.
    \label{fig:FlowChart}
\end{figure}
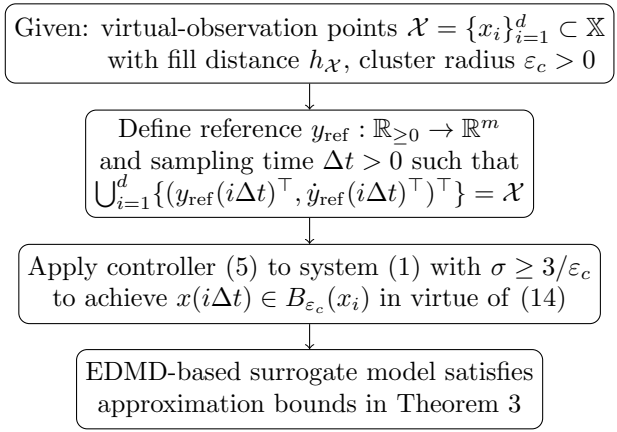

Next, we sketch a scheme, in which the safeguarding controller~\eqref{eq:FunnelControl} from \Cref{Sec:FunnelControl} is utilized to ensure a sufficiently-small fill distance~$h_{\mathcal{X}}$ and cluster radius~$\varepsilon_c$ and, thus, a desired guaranteed approximation accuracy by suitably collecting data triplets, cf. the flowchart in~\Cref{fig:FlowChart}.
First note that the application of the feedback controller~\eqref{eq:FunnelControl} guarantees~\eqref{eq:ControlObjective}, i.e., output tracking within prescribed bounds on the error $y - y_{\rm ref}$. 
To ensure a given fill distance within a subset of the state space~$\mathbb{X} \subset \R^{2m}$, we make the following observation:
For $x(0) = (y_{\rm ref}(0), \dot y_{\rm ref}(0))^\top$ and constant error margin~$\sigma  > 0$, \cite[Lem.~2.1]{lanza2024_sampleddata} yields for all ~$t \ge 0$
\begin{align*}
    \| e_1(t) \| & \le \frac{\sqrt{5}-1}{2} =: \delta <1 \  \text{and} \ \|e_2(t)\| \le 1.
\end{align*}
Rewriting the auxiliary variable~$e_2$ in~\eqref{eq:FunnelControl}, we obtain
\begin{align*}
    \| \dot e(t) \| = \frac{1}{\sigma} \left\|e_2(t) - \tfrac{e_1(t)}{1-\|e_1(t)\|^2} \right\|  \le \frac{1 + \tfrac{\delta}{1-\delta^2}}{\sigma} = \frac{2}{\sigma},
\end{align*}
where~$0<\sigma \in \R_{> 0}$ can be chosen (arbitrarily) large.
Therefore, invoking $\|y(t) - y_{\rm ref}(t)\| < 1/\sigma$, we can estimate the error between the state~$x(t) \in \R^{2m}$ and the reference and its derivative by
\begin{equation} \label{eq:FunnelBoundariesStatespace}
\left\| \begin{pmatrix}
    x_1(t) \\ x_2(t)
\end{pmatrix} 
-
\begin{pmatrix}
    y_{\rm ref}(t) \\ \dot y_{\rm ref}(t)
\end{pmatrix}
\right\| < \frac{3}{\sigma}.
\end{equation}
Based on this line of reasoning, we propose the scheme summarized in \Cref{fig:FlowChart} to generate an EDMD-based surrogate model with prescribed approximation error, where the data is collected by steering the system safely through the state space.
Note that \Cref{fig:FlowChart} is a schematic presentation.
In general, the reference will be designed such that for some $D > d$ we have
$ \bigcup_{i=0}^D\{ (
    y_{\rm ref}(i \Delta t)^\top, \dot y_{\rm ref}(i \Delta t)^\top)^\top \} \supset \mathcal{X}$.
Then, system data is recorded when $x(t) \in B_{\varepsilon_c}$.
Moreover, the error margin~$\sigma$ should be chosen such that multiple data points per cluster can be recorded (cf.~data requirements after~\eqref{eq:DiscreteSystem}).
These aspects are topics of future research.
Provided that the details of the above reasoning are formulated precisely, it can then rigorously be ensured that the EDMD-based MPC controller maintains the desired output-tracking specifications using, e.g., the controller design proposed in~\cite{soloperto2022nonlinear}, meaning that the safeguarding mechanism may be switched off.

\section{Numerical example} \label{Sec:Numerics}
To illustrate the previously presented findings we present a numerical simulation. 
Consider the nonlinear system
\begin{subequations} \label{eq:VdP}
\begin{equation} 
    \begin{aligned}
        \dot x_1(t) &= x_2(t), \\
        \dot x_2(t) &= \nu ( 1- x_1(t)^2)x_2(t) - x_1(t)  +  u(t), 
    \end{aligned}
    \end{equation}
    with output
    \begin{equation}
        y(t) = h(x(t)) = x_1(t),
    \end{equation}
\end{subequations}
where~$\nu = 0.1$, and $(x_1(0),x_2(0))^\top = (1,-1)^\top \in \R^2 $.
System~\eqref{eq:VdP} is a forced Van-der-Pol oscillator. 
In this particular example, we have state dimension~$n=2$ and input-output dimension~$m=1$.
System~\eqref{eq:VdP} satisfies the structural \Cref{Ass:G_PosDef}, and hence belongs to the system class specified in \Cref{subsec:system-class:control-objective}.
According to \Cref{Ass:SignalsAvailable} we have access to the signals~$y(t)$ and~$\dot y(t)$.
We choose the following set of observables
\begin{equation*}
    \Psi = \{ x_1^p \cdot x_2^{q} \ | \ p,q \in \{0,1,2,3\}, \ p+q \le 3 \},
\end{equation*}
meaning that we have~$10$ observables consisting of all monomials of the state up to degree~$3$. In particular, the coordinate functions are contained.
For the cost function in \Cref{Alg:MPC} we choose the weighting matrices $Q= \text{diag}(10^4,1) \in \R^{2 \times 2}$ and $R=10^{-4}$.
We choose the prediction horizon to be~$N = 30$, and constrain the prediction-based input by~$\|\mu(\cdot)\| \le 2$, i.e., $\mathbb{U} = [-2,2]$.
The sampling time for EDMD is set to ~$\Delta t = 0.05$.    
For the activation function~$a_\tau(\cdot)$ in~\eqref{eq:ActivationFunction} we choose~$\lambda = 0.75$ as activation threshold, and $\tau = \Delta t/2$ as dwell-time.

\begin{figure}[htb]
    \centering
    \includegraphics[width=0.9\linewidth]{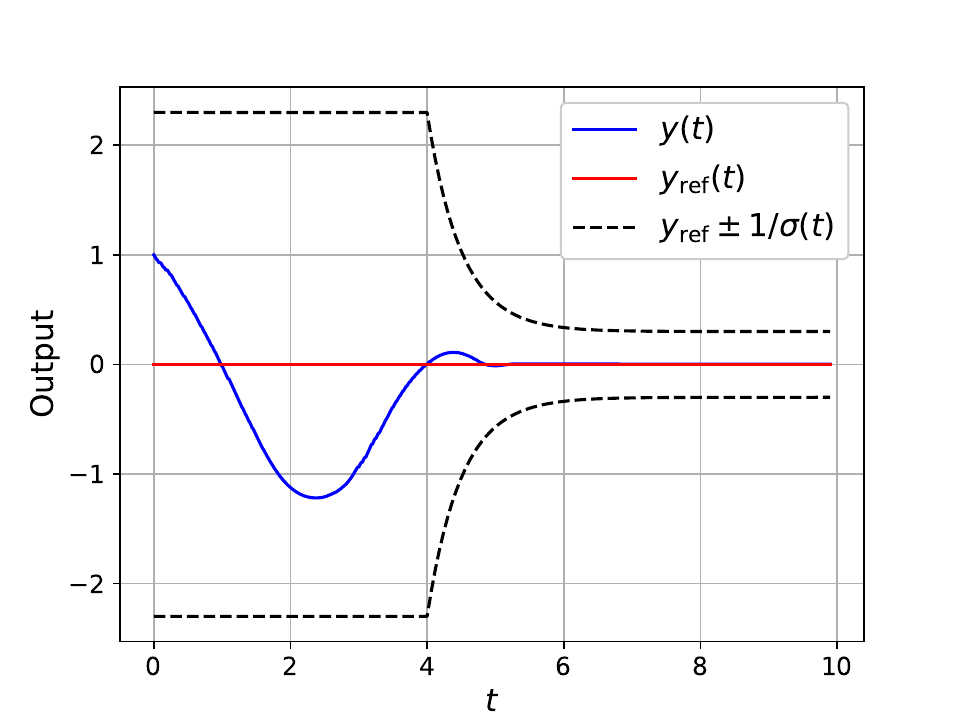}
    \caption{Stabilizing the origin. Visualization of the output~$y(t) = x_1(t)$ and the error boundary~$\sigma(t)$ given by~\eqref{eq:example:error-boundary}.}
    \label{fig:Output_stabilization}
\end{figure}

We consider two control tasks. First we aim to stabilize the system to the origin; second, we perform a set-point transition.
For both tasks we choose the funnel function
\begin{equation}\label{eq:example:error-boundary}
\sigma(t) = \begin{cases}
\frac{1}{2.3} & t \le 4, \\ 
\frac{1}{2e^{-2(t-4)} + 0.3} & t > 4,
\end{cases} 
\end{equation}
which allows a larger error at the beginning, cf.~\Cref{fig:Output_stabilization,fig:Output_trans}.
This prevents early intervention of the safeguard and allows for data collection during tracking.

\subsection{Stabilization}
\begin{figure}[htb]
    \centering
    \includegraphics[width=0.475\linewidth]{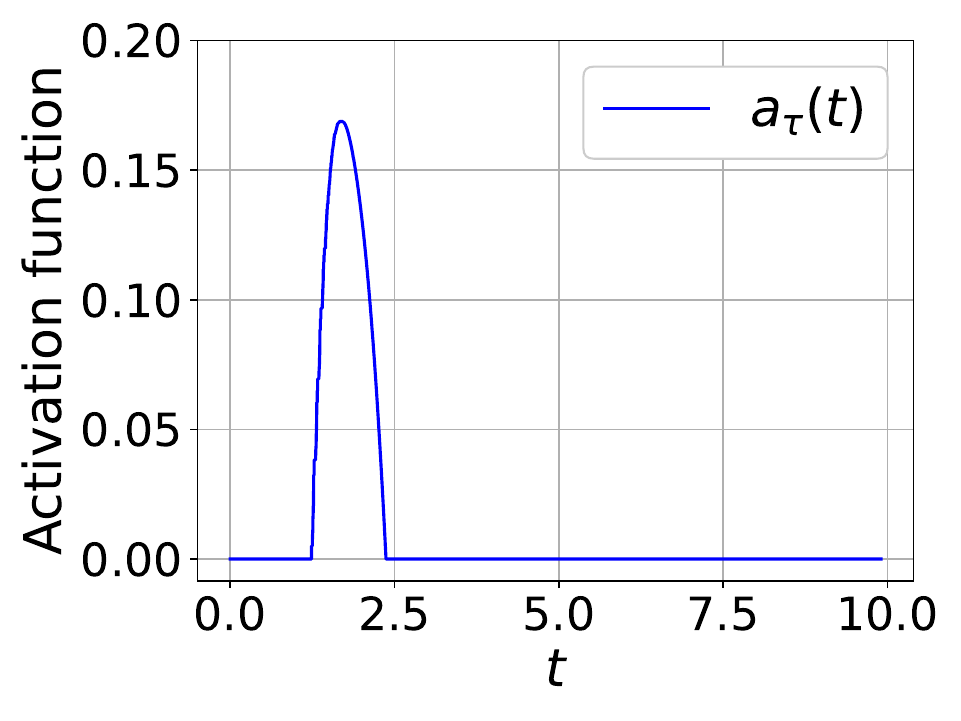}
    \includegraphics[width=0.475\linewidth]{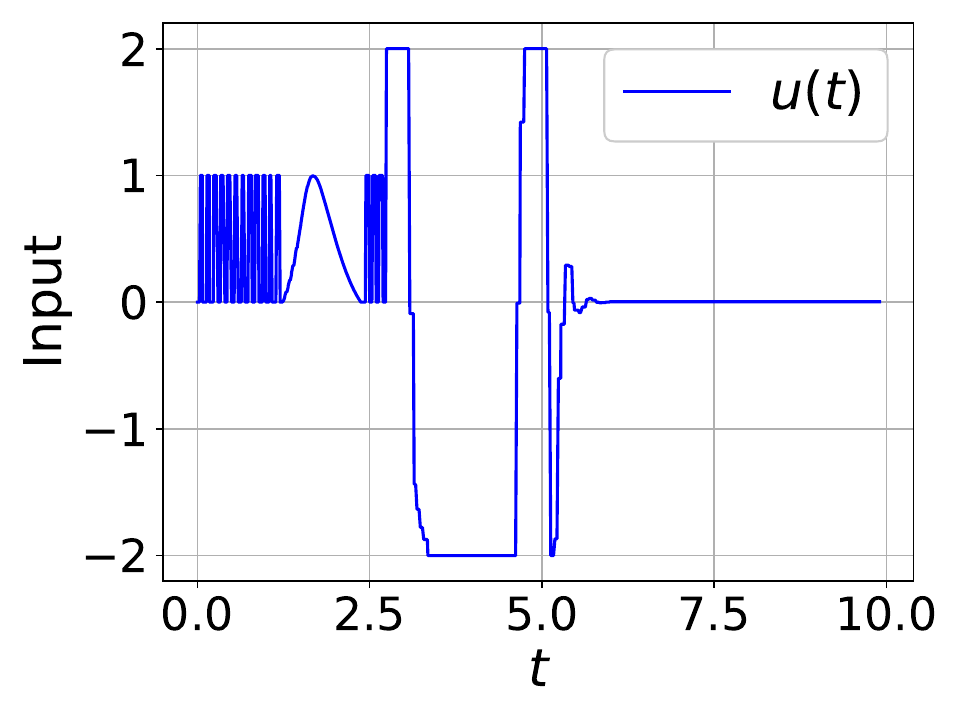}
    \caption{Stabilizing the origin. Visualization of the activation function~$a_\tau(t,e_2(t))$ (left) and the input~$u(t) = \mu(t) + a_\tau(t,e_2(t)) \, u_{\rm FC}(t)$ (right).}
    \label{fig:Control_stabilization}
\end{figure}

First, we simulate stabilization of~\eqref{eq:VdP}, i.e., $y_{\rm ref}(t) \equiv 0$ for all~$t \ge 0$.
We consider the following situation.
We assume that we have access to~$d=10$ data points (e.g. from previous experiments)  collected in~$\mathbb{X} = [2,2]^2$ with sampling time~$\Delta t = 0.05$, i.e., we start with \textit{some} approximation of the Koopman operator.
To improve the model, we collect data whenever $a_\tau(t,e_2(t)) = 0$ for all~$t \in [k \Delta t, (k+1)\Delta t)$ and alternately apply $\mu(t) = 1$ and $\mu(t) = 0$, until~$d = 25$ (cf.~\Cref{rem:aTauZero}).
The control input~$u(\cdot)$ and activation function~$a_\tau(\cdot)$ are depicted in \Cref{fig:Control_stabilization}. 
It can be seen that at some point in the exploration phase the safeguarding funnel controller intervenes.
After some time the approximation of the Koopman operator sufficiently well represents the dynamics and the tracking task is performed with purely data-based predictions.

\subsection{Set-point transition}
As a second scenario we consider transition between two set-points. 
To ensure $y_{\rm ref} \in W^{2,\infty}([0,\infty),\R^2)$, we connect the two-setpoints with a smooth function. 
We define the overall reference trajectory by
\begin{equation*}
    y_{\text{ref}}^{\bar t}(t) = 1 + \tfrac{2}{\pi} \int_0^t e^{-(s-\bar t)^2} \, \text{d}s,
\end{equation*}
parameterized by the shift~$\hat t$.
Note that $y_{\rm ref} \approx 0$ at the beginning, and ${y_{\rm ref} \approx 2}$ in the end.

\subsubsection{Initializing with ten i.i.d. data points}
\begin{figure}[H]
    \centering
    \includegraphics[width=0.475\linewidth]{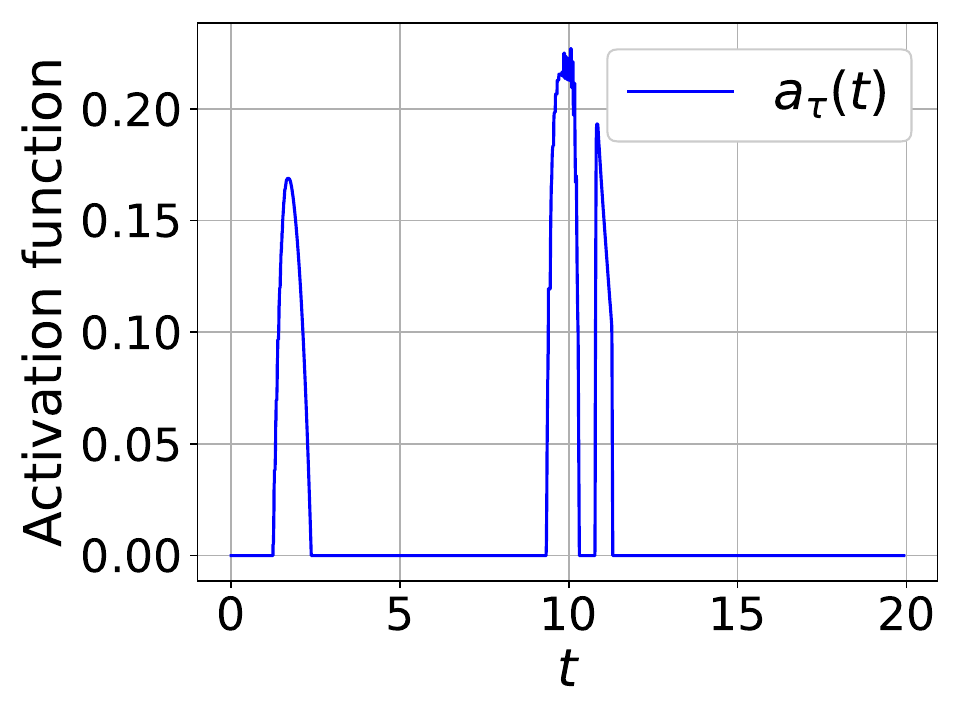}
    \includegraphics[width=0.475\linewidth]{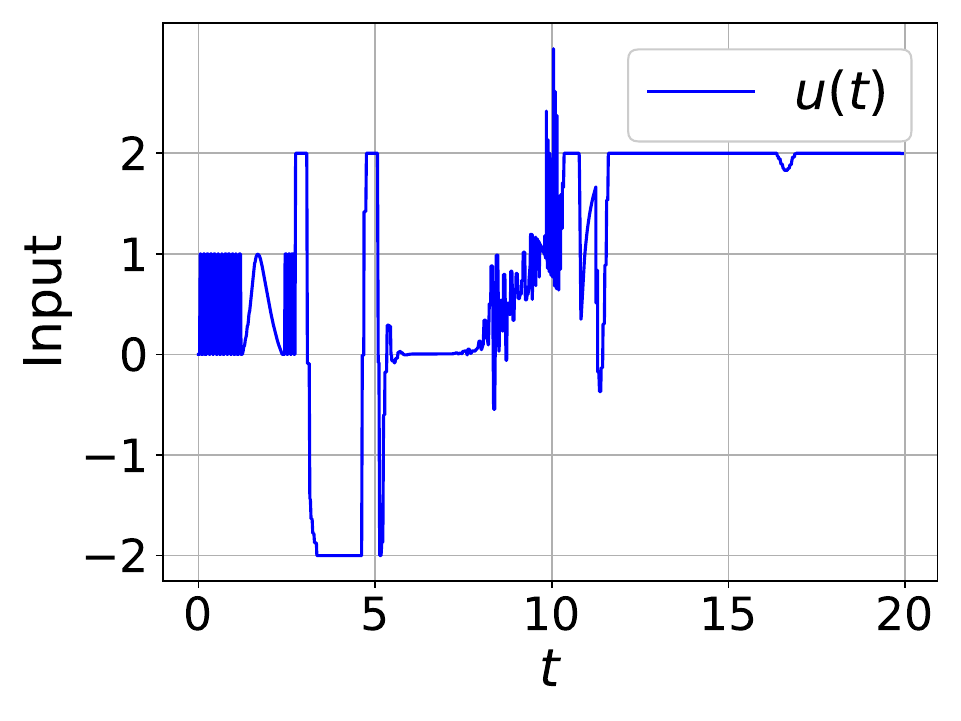}
    \caption{Set-point transition. Visualization of the activation function $a_\tau(t,e_2(t))$ (left) and control input. 
    }
    \label{fig:aTau_trans}
\end{figure}

First we consider the situation when the model is initialized with~$10$ i.i.d. generated data points. 
We choose the reference~$y_{\rm ref}^{10}$.
Accordingly, the results look similar to the previous simulation at the beginning.
\begin{figure}[H]
    \centering
    \includegraphics[width=0.9\linewidth]{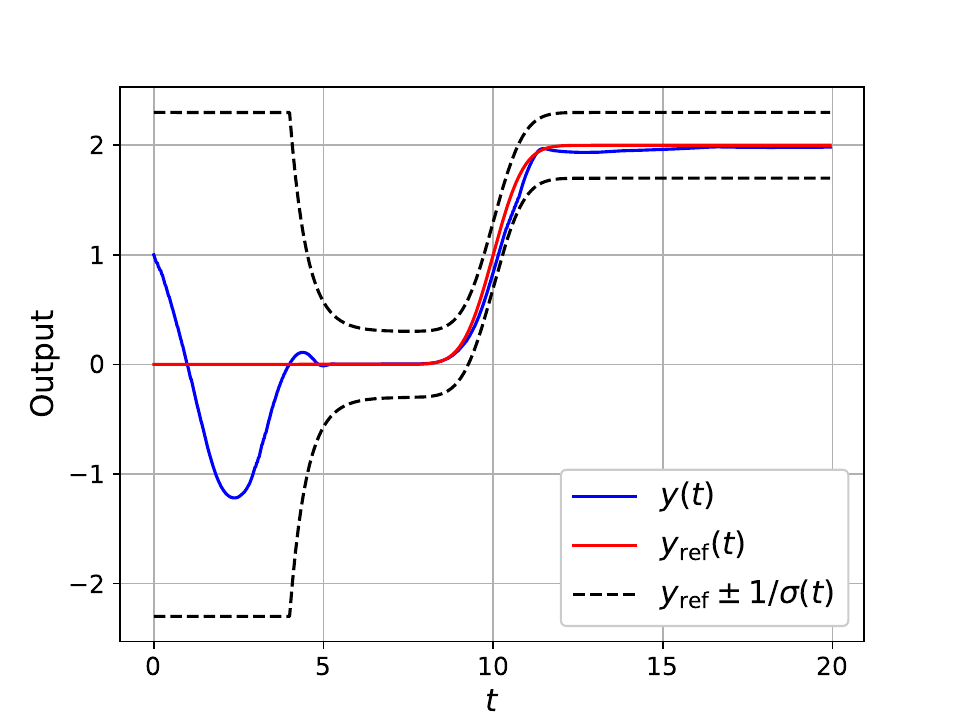}
    \caption{Set-point transition. {Visualization of the} output $y(t) = x_1(t)$, reference and error boundary $1/{\sigma}(t)$.}
     \label{fig:Output_trans}
\end{figure}
In the transition phase maintaining the output constraints is challenging, see \Cref{fig:aTau_trans}, where the funnel controller is activated three times ($a_\tau \neq 0$) and intervenes the predictive control. 
However, the error bounds are guaranteed, see \Cref{fig:Output_trans}.
After the transition phase the EDMD-based MPC is again capable of keeping the system close to the set-point, which is a controlled equilibrium matching the input constraints.

\subsubsection{Initializing with one data point}
\begin{figure}[htb]
    \centering
    \includegraphics[width=0.9\linewidth]{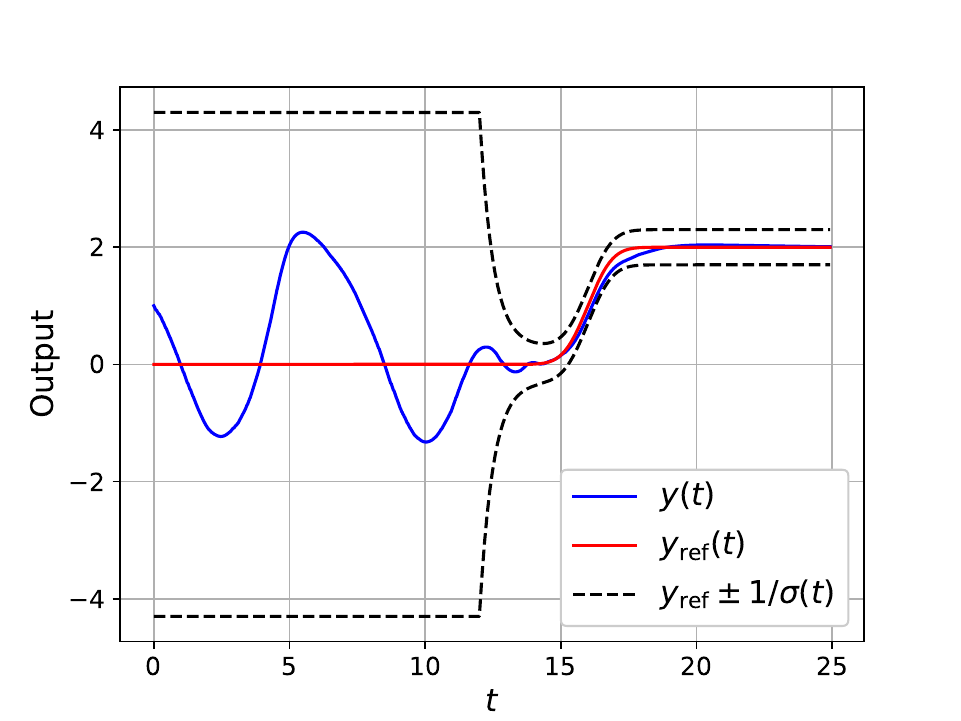}
    \caption{Set-point transition, initializing with 1 data point. Visualization of the output~$y(t) = x_1(t)$, reference and error boundary~$1/\sigma(t)$.}
    \label{Fig:OneData_output}
\end{figure}

\begin{figure}[htb]
    \centering
    \includegraphics[width=0.475\linewidth]{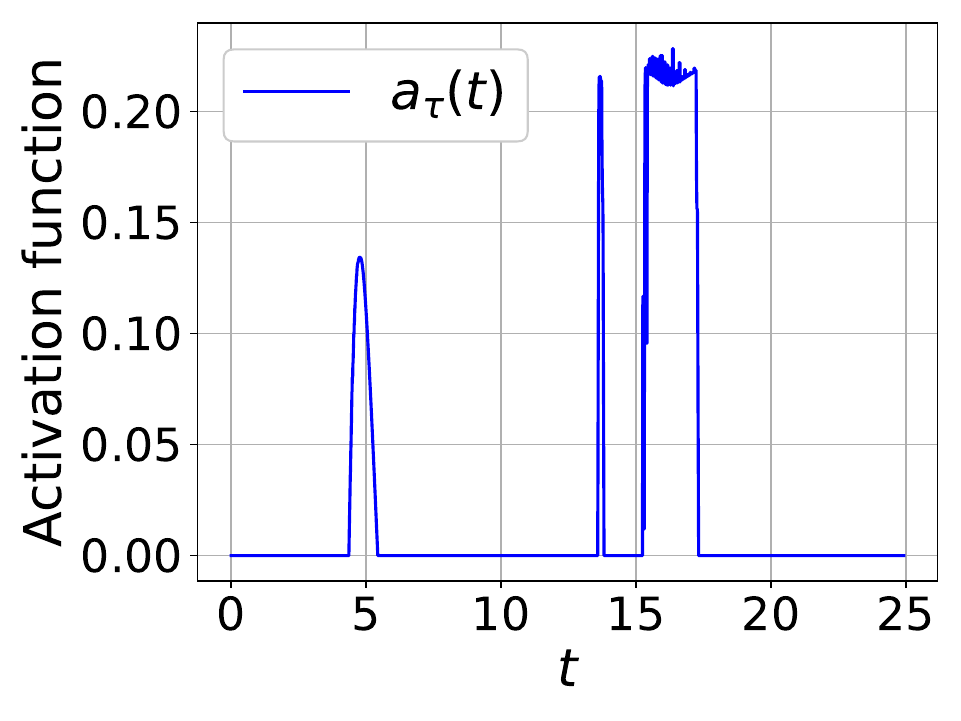}
    \includegraphics[width=0.475\linewidth]{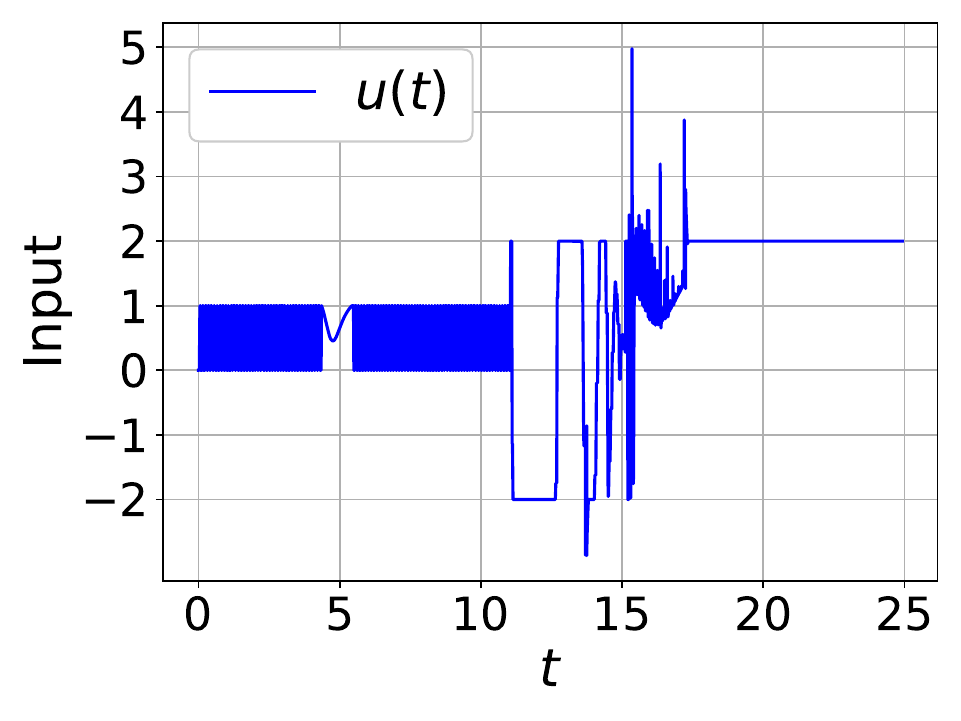}
    \caption{Set-point transition ,initializing with 1 data point. Visualization of the activation function $a_\tau(t,e_2(t))$ (left) and control input (right). 
    }
    \label{Fig:OneData_aTau_trans}
\end{figure}

Next, we start with only one data point and collect data on runtime (the one data point is used to initialize the algorithm).
To account for the lack of data, we start with a larger funnel and choose the reference~$y_{\rm ref}^{16}$, i.e., we allow the system to evolve ``freely'' for a longer period, see \Cref{Fig:OneData_output}.
For simulation purpose we limit the data to~$d = 100$.
Similarly to the previous scenario, the funnel controller intervenes three times in the transition phase due to the tight error margin, cf.~\Cref{Fig:OneData_aTau_trans}.
After the learning process the data-based MPC is capable of keeping the system at the controlled equilibrium.
However, while the error bounds are still maintained during the set-point change, using this data sampling process, leads to control values, that distinctly exceed the input constraints for the MPC.

\section{Conclusions and outlook} 
We have proposed a two-component controller design for nonlinear control systems to achieve safe data-driven predictive tracking control.
The data-driven component has been exemplified by DeePC and EDMD.
Two scenarios have been considered.
First, the reactive model-free feedback controller safeguards the tracking with DeePC and EDMD-based MPC. 
Second, using the guarantees from the feedback can be utilized to achieve a certain sampling quality (fill distance).
The proposed two-component controller~\eqref{eq:ControlSignal} is open to many extensions and modifications such as
updating the data matrices and deleting old data (shifting the sampling window);
starting the control task with a small prediction horizon and depending on the amount of data, the prediction horizon is gradually increased (as done for DeePC in~\cite{schmitz2023safe}); or use only few data (cf. model order reduction), which results in small data matrices and hence saves computation time.
Depending on the application, these ideas can be selected and combined while maintaining the guarantees for the tracking error.
In future work we will thoroughly analyze the scheme proposed in \Cref{Sec:PrescribedApproximationError} to, e.g., formulate requirements on the data and the reference defined through the virtual-observation points.

\textbf{Acknowledgment.}
    L.\ Bold and L.\ Lanza gratefully acknowledge funding by the German Research Foundation DFG (Project-ID 471539468 and 545246093; Deutsche Forschungsgemeinschaft) and by the Carl Zeiss Foundation (Project-ID~2011640173; VerneDCt).

\printbibliography
\end{document}